\documentclass[12pt,oneside,english]{amsart}
\usepackage{bbding}
\usepackage{dsfont}
\usepackage{graphicx}
\usepackage{amsthm,amscd}
\usepackage{calc}
\usepackage{enumerate}
\usepackage{mathrsfs}
\usepackage{mathrsfs}
\usepackage{CJK,fancyhdr,amscd}
\usepackage{anysize}
\usepackage{amsmath,amssymb,amsfonts}
\usepackage{epsfig}
\usepackage{latexsym}
\usepackage[pagebackref,colorlinks]{hyperref}
\usepackage{indentfirst, latexsym}
\usepackage{graphics}
\usepackage{amsmath, bm}
\usepackage{amsfonts}
\usepackage{amssymb}%
\usepackage{color}
\usepackage{mathtools}
\usepackage{verbatim}

\providecommand{\U}[1]{\protect\rule{.1in}{.1in}}

\numberwithin{equation}{section}
\providecommand{\U}[1]{\protect\rule{.1in}{.1in}}

\newtheorem{theorem} {Theorem} [section]
\newtheorem{proposition}[theorem]{Proposition}
\newtheorem{corollary}  [theorem]     {Corollary}
\newtheorem{lemma}  [theorem]     {Lemma}
\newtheorem{example}  [theorem]     {Example}
\newtheorem{remark}  [theorem]     {Remark}
\newtheorem{definition}  [theorem]     {Definition}

\newtheorem{question}  [theorem]     {Question}

\newcommand{\bproof}{\begin{proof}}
\newcommand{\eproof}{\end{proof}}

\newcommand{\p}{\partial}
\newcommand{\bp}{\bar\partial}

\renewcommand{\qedsymbol}{$\square$}
\newcommand{\bthm}{\begin{theorem}}
\newcommand{\ethm}{\end{theorem}}
\newcommand{\blem}{\begin{lemma}}
\newcommand{\elem}{\end{lemma}}
\newcommand{\bprop}{\begin{proposition}}
\newcommand{\eprop}{\end{proposition}}
\newcommand{\brmk}{\begin{remark}}
\newcommand{\ermk}{\end{remark}}
\newcommand{\bcor}{\begin{corollary}}
\newcommand{\ecor}{\end{corollary}}
\newcommand{\bdefi}{\begin{definition}}
\newcommand{\edefi}{\end{definition}}
\newcommand{\beq}{\begin{equation}}
\newcommand{\eeq}{\end{equation}}
\newcommand{\bex}{\begin{example}}
\newcommand{\eex}{\end{example}}

\newcommand{\ov}{\overline}
\newcommand{\ti}{\tilde}
\usepackage[lmargin=1in, rmargin=1in]{geometry}
\title{The twisted constant in Calabi-Yau type equation   }
\author{Genglong Lin}
\address{Bejing Institute of Mathematical Sciences and Application, Beijing 101408, China}
\email{lingenglong@amss.ac.cn}

\begin{document}
\date{\today}
\begin{abstract}
In this paper we establish a necessary and sufficient condition for solving a 
general class of fully nonlinear elliptic equations on compact almost Hermitian manifolds,
extending a recent work of Guo-Song. As applications, we determine the twisted constants 
in Calabi-Yau type equations, including the classical one, Hessian type one and form type one with gradient terms introduced by Popovici and Tosatti-Weinkove.
In particular, it addresses a question raised in a work of 
Chu-Tosatti-Weinkove \cite[Introduction, Remark 5]{CTW19}. 
\end{abstract}
\maketitle

\section{introduction}\label{introduction}
Let $(M,\chi,J)$ be a compact almost Hermitian manifold of real dimension
$2n$ and $\omega$ be a fixed real $(1,1)$ form on $(M,J)$. For a smooth
function $u$ on $M$, we write 
$$\omega_u:=\omega+\sqrt{-1}\p\bp u+Z(\p u)=\omega+\frac{1}{2}(dJdu)^{1,1}
+Z(\p u),$$
where $Z(\p u)$ is a smooth $(1,1)$ form which depends on $\p u$ linearly.
Let $\lambda(u)=(\lambda_1(u),...,\lambda_n(u))$ be the eigenvalues of $\omega_u$ with 
respect to $\chi$. We consider a general class of PDE as follows:
\begin{equation}\label{eq1.1}
    F[u]:=f(\lambda_1,...,\lambda_n)=e^{\psi},
\end{equation}
where $\psi\in C^{\infty}(M)$ and $f$ is a 
smooth symmetric function defined on an open convex cone $\Gamma\subset\mathbb{R}^n$
with vertex at the origin. We also assume that $f\in C^{\infty}(\Gamma)\cap C^0(\ov{\Gamma})$
and the positive cone $\Gamma^n\subset\Gamma$. Furthermore, $f$ satisfies the 
following structural conditions:
\begin{enumerate}[a)]
    \item \label{1.2} $\frac{\p f}{\p \lambda_i}>0$ in $\Gamma,i=1,...,n,$\\
    \item\label{1.3} $f$ is concave in $\Gamma$,\\
    \item\label{1.4} $f>0$ in $\Gamma$, $f=0$ in $\p\Gamma$,\\
    \item\label{1.5} $\lim_{R\to\infty}f(R\lambda)=\infty$ for any $\lambda\in\Gamma$.
\end{enumerate}
There are many examples of $f$ satisfying the conditions above such as 
$$ f(\lambda)=(S_k(\lambda))^{1/k}=\left(\sum_{i_1<i_2<...<i_k}\lambda_{i_1}\lambda_{i_2}...\lambda_{i_k}\right)^{1/k},$$
for $k=1,...,n$. The corresponding cone $\Gamma$ is given by 
$$\Gamma_k=\{\lambda\in\mathbb{R}^n|S_i(\lambda)>0,i=1,...,k\}.$$
When additionally $Z=0,k=n$, $(M,\chi,J)$ is integrable and $\omega$ is a K\"ahler form, this 
corresponds to the classical Calabi-Yau equation which can be solved since Yau's
solution to Calabi's conjecture \cite{Yau78}. The form type Calabi-Yau equation for $(n-1)$-plurisubharmonic
functions was introduced and studied by Fu-Wang-Wu \cite{FWW10}\cite{FWW15}, which can be 
formulated as follows:
\begin{equation}\label{form type eq}
    \left(\eta+\frac{1}{n-1}\left((\Delta^C u)\chi-\sqrt{-1}\p\bp u\right)\right)^n=e^{\psi}\chi^n,
\end{equation}
where $\eta$ is a Hermitian metric and $\Delta^C$ denotes the canonical Laplacian operator of $\chi$.
When $\chi$ is K\"ahler, the equation (\ref{form type eq})  was solved by Tosatti-Weinkove \cite{TW10b} 
and the result was generalized to general Hermitian metric $\chi$ in \cite{TW19}. A more general 
type equation, which contained gradient terms of $u$ as follows, was introduced by Popovici \cite{Popovici15}
and Tosatti-Weinkove \cite{TW19} and solved by Szekelyhidi-Tosatti-Weinkove \cite{STW17}, confirming a conjecture of Gauduchon:
\begin{equation}\label{form type eq1}
    \left(\eta+\frac{1}{n-1}\left((\Delta^C u)\chi-\sqrt{-1}\p\bp u\right)+W(\p u)\right)^n=e^{\psi}\chi^n.
\end{equation}
Here $W$ is a Hermitian tensor that linearly depends on $\p u$. Also see \cite{TW17,TW21}.

The complex Hessian case $1\leq k\leq n,Z=0$ was solved by Dinew-Kolodziej \cite{DK17} on compact K\"ahler manifold and was 
then generalized to general Hermitian manifold by Szekelyhidi \cite{Szekelyhidi18} and Zhang \cite{Zhang17}
independently. 

Calabi-Yau type equations on non-integrable almost hermitian 
manifolds have also aroused the interest of many people in the subject over the last decade. 
A remarkable result in this direction is due to Chu-Tosatti-Weinkove \cite{CTW19}.
They proved the following theorem:
\begin{theorem}\cite{CTW19}
    Let $(M,\chi,J)$ be a compact almost Hermitian manifold of real dimension $2n$. Given a smooth 
    funtion $h$ there exists a unqiue pair $(u,b)$ where $u\in C^{\infty}(M)$ and $b\in\mathbb{R}^n$,
    solving 
    \begin{equation}
        \begin{aligned}
            &(\chi+\sqrt{-1}\p\bp u)^n={e^{\psi+c}}\chi^n,\\
& \chi+\sqrt{-1}\p\bp u>0,\sup_M u=0.
        \end{aligned}
\end{equation}
\end{theorem}
An interesting question is to determine the constant $b$, which we call 
\emph{the twisted constant}, from $\chi$ and  $h$ via a simple formula. 
An easy application of maximum principle implies that $|b|\leq\sup_M|F|$,
as noted by Chu-Tosatti-Weinkove \cite[Introduction, Remark 5]{CTW19}.
 
The results on solutions of Hessian type equation and form type Calabi-Yau equation 
with gradient terms on compact Hermitian manifolds were extended to 
almost Hermitian case by Chu-Huang-Zhang \cite{CHZ23} and Huang-Zhang \cite{HZ23}.
It is also interesting to ask:
\begin{question}
    Let $(M,\chi,J)$ be a compact almost Hermitian manifold of real dimension $2n$.
    \begin{enumerate}
        \item Assume that $\eta$ is an almost Hermitian metric and $(u,b)\in C^{\infty}(M)\times\mathbb{R}$ such 
        that
       \begin{equation}
        \begin{cases}
            \left(\eta+\frac{1}{n-1}\left((\Delta^C u)\chi-\sqrt{-1}\p\bp u\right)+W(\p u)\right)^n={e^{\psi+c}}\chi^n,\\
            \left(\eta+\frac{1}{n-1}\left((\Delta^C u)\chi-\sqrt{-1}\p\bp u\right)+W(\p u)\right)>0,\\
            \sup_M u=0.
        \end{cases}
       \end{equation}
       \item Assume that $\omega$ is a smooth $k$-positive real $(1,1)$-form and 
       $1\leq k\leq n$. Let $(u,b)\in C^{\infty}(M)\times\mathbb{R}$ be such that 
       \begin{equation}
        \begin{cases}
            \omega_u^k\wedge\chi^{n-k}={e^{\psi+c}}\chi^n,\\
            \frac{\omega_u^i\wedge\chi^{n-i}}{\chi^n}>0,i=1,2,...,k,\\
            \sup_M u=0.
        \end{cases}
       \end{equation}
    \end{enumerate}
Determine the twisted constants $c$ above.
\end{question}
Recently, Guo-Song \cite{GS24} established an analytic criterion for the 
solvablity problem of (\ref{eq1.1}) where the gradient terms vanished
 on compact Hermitian manifold and determined 
$b$ in a simple way. The aim of this paper is to extend their result to a more general class
of equations which may contain gradient terms on compact almost Hermitian 
manifolds and give a simple expression for $c$. Namely, we can answer the question above as following:
\begin{theorem}\label{mainthm}
    In case $(1)$, $$e^c=\inf_u \max_M \frac{e^{-\psi}(\eta+\frac{1}{n-1}\left((\Delta^C u)\chi-\sqrt{-1}\p\bp u\right)+W(\p u))^n}{\chi^n},$$
    where the infimum is taken on the set $\left\{u\in C^{\infty}(M),\eta+\frac{1}{n-1}\left( (\Delta^C u)\chi-\sqrt{-1}\p\bp u\right)+W(\p u)>0\right\}$.\\
            In case $(2)$, $$e^c=\inf_u\max_M \frac{e^{-\psi}\omega_u^k\wedge\chi^{n-k}}{\chi^n},$$ where 
            the infimum is taken on the set
            $\left\{u\in C^{\infty}(M),\frac{\omega_u^i\wedge\chi^{n-i}}{\chi^n}>0,i=1,2,...,k\right\}.$
\end{theorem} 
The theorem above is based on a similar criterion to Guo-Song's on solvablity of (\ref{eq1.1}). We remark 
that compared to Hermitian setting, the proof requires more delicate analysis and overcoming some difficulties
which come from the non-integrable almost complex structure on the ambient space and gradient terms in the equation. 
  
Our paper is organized as follows. In Section \ref{sec1} we will introduce some basic notions including almost complex 
structures, sub-solutions and sub-slopes. In Section \ref{sec2} we will 
prove a Guo-Song type criterion on solvablity of a class of fully nonlinear elliptic equations (\ref{eq1.1})
which contain gradient terms, and then complete the proof of theorem above.

\section{preliminary}\label{sec1}
Suppose that $(M,\chi,J)$ is an almost Hermitian manifold of real dimension $2n$.
From basic knowledge of almost Hermitian manifolds in the literature, one can define by duality
$(p,q)$-forms and operator $\p,\bp$ using the decomposition of complexified tangent
space $T^{\mathbb{C}}M$ into the two eigenspaces $T^{1,0}M, T^{0,1}M$.
Let $A^{1,1}(M)$ denote the set of smooth real $(1,1)$ forms on $(M,J)$.
For any $u\in C^{\infty}(M)$, we have that $\sqrt{-1}\p\bp u=\frac{1}{2}(dJdu)^{1,1}$
is a real $(1,1)$ form in $A^{1,1}(M)$. Let $\omega$ be a real $(1,1)$ form on $M$
and set $$\omega_u=\omega+\sqrt{-1}\p\bp u+Z(\p u),$$
where $\p u$ is a real $(1,1)$ form defined by $Z_{i,\bar{j}}=Z_{i,\bar{j}}^pu_p+\overline{Z_{i\bar{j}}^p}u_{\bar{p}}$
in a local frame.

For any point $x_0\in M$, let $(e_1,...,e_n)$ be a local unitary $(1,0)$ frame with
respect to the almost hermitian metric $\chi$ near $x_0$ and $\{\theta^i\}_i$ be
the dual coframe. Then locally we can write
$$\chi=\sqrt{-1}\delta_{i\bar{j}}\theta^i\wedge\overline{\theta}^j$$
and $$\omega=\sqrt{-1}g_{i\bar{j}}\theta^i\wedge\overline{\theta}^j,
\omega_u=\sqrt{-1}\tilde{g}_{i\bar{j}}\theta^i\wedge\overline{\theta}^j,$$
where $$\tilde{g}_{i\bar{j}}=g_{i\bar{j}}+e_i\overline{e}_j (u)-[e_i,\overline{e}_j]^{(0,1)}(u)+u_p Z_{i\bar{j}}^p+u_{\bar{p}}\ov{Z_{i\bar{j}}^p},$$
and $[e_i,\ov  {e}_j]^{(0,1)}$ is the $(0,1)$ part of the Lie bracket $[e_i,\ov  {e}_j]$.

As in the literature, we define 
\begin{equation*}
    G^{i\overline{j}}=\frac{\partial F}{\p\tilde{g}_{i\ov  {j}}}, \qquad
    G^{i\overline{j},k\overline{l}}=\frac{\partial^{2}F}{\p\ti{g}_{i\ov  {j}}\ti{g}_{k\ov  {l}}}.
\end{equation*}
After making a unitary transformation, we may assume that  $\tilde{g}_{i\overline{j}}(x_{0})=\delta_{ij}\tilde{g}_{i\overline{i}}(x_{0})$.
We denote $\tilde{g}_{i\overline{i}}(x_{0})$ by $\mu_{i}$.
It is useful to order $\mu_{i}$ such that
\begin{equation}\label{mu order}
    \mu_{1}\geq\mu_{2}\geq\cdots\geq\mu_{n}.
\end{equation}
At $x_{0}$, we have the expressions of $G^{i\ov{j}}$ and $G^{i\bar{k},j\bar{l}}$ (see e.g. \cite{Andrews94,Gerhardt96,Spruck05})
	\begin{equation}\label{second derive of F}
		G^{i\ov{j}} = \delta_{ij}f_{i},\qquad
		G^{i\bar{k},j\bar{l}}=f_{ij}\delta_{ik}\delta_{jl}
		+\frac{f_{i}-f_{k}}{\mu_{i}-\mu_{k}}
		(1-\delta_{ik})\delta_{il}\delta_{jk},
	\end{equation}
	where the quotient is interpreted as a limit if $\mu_{i}=\mu_{j}$. Using \eqref{mu order}, we obtain (see e.g. \cite{EH89,Spruck05})
	\[
	G^{1\ov{1}} \leq G^{2\ov{2}} \leq \cdots \leq G^{n\ov{n}}.
	\]
	On the other hand, the linearized operator of equation \eqref{eq1.1} is
	\begin{equation}\label{L}
		L(v)=G^{i\bar{j}}\Big(e_{i}\bar{e}_{j}(v)-[e_{i},\bar{e}_{j}]^{0,1}(v)+e_{p}(v)Z_{i\bar j}^{p}+\bar{e}_{p}(v)\overline{Z_{i\bar j}^{p}}\Big).
	\end{equation}
    Since $[e_{i},\bar{e}_{j}]^{(0,1)}$ is a differential operator of first order, $L$ is a second order elliptic operator.
\textbf{Sub-slopes and sub-solutions}\\
As in \cite{GS24}, we define $\mathcal{E}$ to be the set of admissible functions 
for equation (\ref{eq1.1}) associated to the cone $\Gamma$ by 
 \begin{equation}\label{eqn:e definition}
 \mathcal{E}  =\mathcal{E}  (M, \chi, \omega, f, \Gamma)= \{ u \in C^\infty(M)~|~ \lambda(\omega_u(p)) \in \Gamma, ~{\rm for~any}~ p\in M \}.
 \end{equation}
We assume that $\mathcal E$ is not empty.  
 \begin{definition}\label{gslope} 

 The sup-slope $\sigma$ for equation (\ref{eq1.1})  associated to the cone $\Gamma$ is defined to be 
 \begin{equation}\label{supslo}
 \sigma = \inf_{ u \in \mathcal{E}} \max_M   e^{-\psi} F[u]   \in [0, \infty). 
 \end{equation} 
 \end{definition}
 In fact, $\sigma$ is always positive if  $\mathcal{E}\neq \emptyset$. Also define $f_{\infty, i}: \Gamma  \rightarrow (0, \infty]$ by  
\begin{equation}
f_{\infty, i} (\lambda_1, ..., \lambda_i, ...,      \lambda_{n}) = \lim_{R\rightarrow \infty} f(\lambda_1, ..., \lambda_{i-1},  R, \lambda_{i+1}, ..., \lambda_n)
\end{equation}
for $i = 1, 2, ..., n$ and 
$$f_\infty (\lambda) = \min_{i=1, ..., n}  f_{\infty, i}(\lambda) .$$
\begin{remark}
    Either $f_\infty \equiv \infty$ on $\Gamma$ or $f_\infty$ is a positive symmetric concave function in $\Gamma$ satisfying (\ref{1.2}) and  (\ref{1.5}).
\end{remark}

\begin{definition}The global subsolution operator 
$$F_\infty:\mathcal{E} \rightarrow C^0(M) \cup \{\infty\}$$ 
is define by
\begin{equation}
F_\infty[u] = f_\infty(\lambda(u)). 
\end{equation}
\end{definition}
It is easy to see $f_{\infty}(\lambda(u))$ is a globally defined function 
on $M$.

Now we can introduce a set of sub-solutions paired with 
the sup-slope for equation (\ref{eq1.1}):
\begin{definition} \label{subdef} Let $\sigma$ be the sup-slope for equation (\ref{eq1.1}) in Definition \ref{gslope}. Then $\underline u \in \mathcal{E} $ is said to be a  sub-solution  associated with  $\sigma$, if 
\begin{equation}\label{subsol}
  e^{-\psi} F_\infty[ \underline u]   > \sigma, 
\end{equation}
on $M$.
\end{definition}
The criterion of solvability can be now stated as follows:
\begin{theorem}\label{criterion}
    Let $(M,\chi,J)$ be a compact almost Hermitian manifold of real dimension $2n$
    and $f$ satisfies the structure conditions (\ref{1.2}), (\ref{1.3}), (\ref{1.4}),
    (\ref{1.5}) for an open convex symmetric cone $\Gamma\subset\mathbb{R}^n$ containing 
    $\Gamma^n$. The the following statement are equivalent:
    \begin{enumerate} 

\item There exists a smooth solution $u \in  $ solving equation (\ref{eq1.1}), or equivalently, 
$$ e^{-\psi} F[u] = constant. $$

\item  There exists a sub-solution $\underline u \in \mathcal{E}$ for equation (\ref{eq1.1}) associated with the sup-slope $\sigma$, i.e., 
$$e^{-\psi} F_\infty[\underline u] > \sigma .$$

\item There exists  $u \in \mathcal{E}$ satisfying
$$ \max_{ M} e^{-\psi}  F[u]  <  \min_{M  } e^{-\psi } F_{\infty}[u].$$

\item There exist  $\overline{u}, \underline u \in \mathcal{E}$ satisfying
$$ \max_M e^{-\psi} F[ \overline u]  <  \min_{M  } e^{-\psi} F_{\infty} [\underline u] .$$

\end{enumerate}
Furthermore, if $u \in \mathcal{E}$ solves equation (\ref{eq1.1}), $u$ is unique up to a constant and   
$$F[u] = \sigma e^\psi,$$
where $\sigma$ is the sup-slope.
\end{theorem}

    \section{Proof of main theorem}\label{sec2}
In this section we use  technique  of Guo-Song \cite{GS24} and continuity method to prove Theorem \ref{criterion}. The 
steps are $(4)\implies (1)\implies (3)\implies (2)\implies (4)$, in which the most
essential one is $(4)\implies (1)$. Now let $\ov{u}$ and $\underline{u}\in\mathcal{E}$ be a pair of super
and sub-solutions satisfying 
\begin{equation}\label{2.1}
    e^{\ov{c}}=\max_M e^{-\psi}f(\lambda (\omega_{\ov{u}}))<\min e^{-\psi}f_{\infty}(\lambda(\omega_{\underline{u}})).
\end{equation}
Also let $\sigma $ be the sup-slope of the equation (1.6), we have $\sigma\leq e^{\ov{c}}$.
It is easy to see that there exists a number $\delta>0$ satisfying
\begin{equation}\label{2.2}
    e^{\ov{c}}<(1+\delta)e^{\ov{c}}\leq\min_M e^{-\psi}f_{\infty}(\lambda(\omega_{\ov{u}})).
\end{equation}
Let $\ov{\psi}=\log F[\ov{u}]$. Then we have $F[\ov{u}]=e^{\ov{\psi}}\leq e^{\psi+\ov{c}}$
by the definition of $\ov{c}$. Immediately we obtain that $\ov{\psi}\leq\psi+\ov{c}$.

Consider the following equations with respect to a parameter $t\in [0,1]$:
\begin{equation}\label{2.4}
    F[\ov{u}+\phi_t]=e^{\psi_t+c_t},\psi_t=(1-t)\ov{\psi}+t\psi,
\end{equation}
where $\ov{u}+\phi_t\in\mathcal{E},\sup_M\phi_t=0$ and $c_t$ is a constant for 
each $t$.

Let $\mathcal{T}=\{t\in [0,1]|(\ref{2.4})~admits~a~smooth~solution~at~t\}.$ Then we 
can prove that 
\begin{lemma}
    $\mathcal{T}$ is non-empty and open.
\end{lemma}
\begin{proof}
    Obviously the equation (\ref{2.4}) can be solved with $\phi_0=0$ and
    $c_0=0$ when $t=0$. Therefore $\mathcal{T}$ is non-empty. Now we prove 
    the openness. Suppose that $\hat{t}\in\mathcal{T}$. Denote the linearized 
    operator of $\log F[\ov{u}+\cdot]$ at $\phi_{\hat{t}}$ by $L_{\phi_{\hat{t}}}(\phi)$ for $\phi\in C^2(M)$.
   It is easy to verify that $\log F$ satisfies the structural condition (\ref{1.2}) and (\ref{1.3}).
   Hence $L_{\phi_{\hat{t}}}$ is a second order elliptic operator and is
    of index zero. Indeed, although $L_{\phi_{\hat{t}}}$ contains some first order terms,
     it is still homotopic to the Laplacian operator and 
    hence has the same index as Laplacian operator, which is zero. This 
    can be easily seen through the Atiyah-Singer index theorem.

    By the maximum principle, $Ker(L_{\phi_{\hat{t}}})=\{constant\}$.
    Denote $L_{\phi_{\hat{t}}}^{*}$ by the $L^2$-adjoint operator of 
    $L_{\phi_{\hat{t}}}$ with respect to the volume form $dV$.
    By the Fredholm alternative and strong maximum principle, 
    $Ker(L_{\phi_{\hat{t}}}^{*})=span \{\xi\}$ where $\xi>0$ is a smooth
    function. Up to a constant, we may assume that $$\int_M \xi dV=1.$$
    
    Define the Banach space by 
    $$\mathcal{U}^{2,\alpha}:=\left\{ g\in C^{2,\alpha}(M)\Big|g\in\mathcal{E},\int_M g\cdot\xi dV=0\right\}.$$
    Its tangent space is given by $$T_{\phi_{\hat{t}}}\mathcal{U}^{2,\alpha}
    :=\left\{h\in C^{2,\alpha}(M)\Big|\int_M h\cdot\xi dV=0\right\}.$$
    Consider the map $$G(\phi,c)=\log F(\ov{u}+\phi_{{t}})- c,$$
    mapping $\mathcal{U}^{2,\alpha}\times\mathbb{R}$ to $C^{\alpha}(M)$.
    The linearized operator of $G$ is $$L_{\phi_{\hat{t}}}- c:T_{\phi_{\hat{t}}}\mathcal{U}^{2,\alpha}\times\mathbb{R}
    \to C^{\alpha}(M).$$ Since for any $w\in C^{\alpha}(M)$ there exists
    a constant $C$ such that $\int_M (w+C)\cdot\xi=0$. By the Fredholm alternative
    again, there exists a real function $\phi$ on $M$ such that 
    $$L_{\phi_{\hat{t}}}(\phi)=w-c.$$
    Hence the map $L_{\phi_{\hat{t}}}- c$ is surjective. An easy application 
    of maximum principle can deduce the injectivity. By the inverse function
    theorem, when $t$ is close to $\hat{t}$, there exists a pair $(\phi_t,c_t)\in\mathcal{U}^{2,\alpha}\times\mathbb{R}$
    such that $$G(\phi_t,c_t)=\log F[\ov{u}+\phi_t]-c=t\psi+(1-t)\ov{\psi}.$$
    The openness is proved.
\end{proof}

Proofs of following lemmas can be found in \cite{GS24}, hence we omit them.
\begin{lemma}\label{lem2.3}
    There exists $C>0$ such that for all $t\in\mathcal{T}$, $$-C\leq c_t\leq t\ov{c}.$$
\end{lemma}
The following lemmas presents basic properties of $f_{\infty}$ and the subsolution operator $F_{\infty}$, 
in which case the concavity is used in the proof.
\begin{lemma}\label{lem3.1}
    Either $f_{\infty}\equiv \infty$ in $\Gamma $ or $f_{\infty}(\lambda)$ is 
    bounded for each $\lambda\in\Gamma$.
\end{lemma}
As a direct consequence of Lemma \ref{lem3.1}, we have 
\begin{corollary}
    Either $F_{\infty}[u]\equiv\infty$ on $M$ for each $u\in\mathcal{E}$
    or $F_{\infty}[u]\in C^{0}(M)$ for each $u\in\mathcal{E}$.
\end{corollary}
\begin{lemma}\label{lem3.4}
    For any $\psi\in C^{\infty}(M)$ and $C\in\mathbb{R}$, the sup-slopes
    satisfies $$\sigma (M,\chi,\omega,f,\psi+C)=e^{-C}\sigma(M,\chi,\omega,f,\psi)>0.$$
    Moreover, for any $\psi_1,\psi_2\in C^{\infty}(M)$, we have 
    $$e^{-\max_M |\psi_2-\psi_1|}\leq\frac{\sigma(M,\chi,\omega,f,\psi_2)}{\sigma(M,\chi,\omega,f,\psi_1)}\leq e^{\max_M|\psi_2-\psi_1|}.$$
\end{lemma}
Using (\ref{2.1}) and Lemma \ref{lem2.3}, we have 
\begin{proposition}\label{3.1}
    Let $\delta>0$ be the fixed constant in (\ref{2.2}) and $\underline{u}\in\mathcal{E}$
    be the sub-solution defined in (\ref{2.2}). Then for any $t\in\mathcal{T}$, 
    we have $$\min_M\left(e^{-\psi_t}f_{\infty}(\lambda(\omega_{\underline{u}}))\right)
    \geq (1+\delta)e^{c_t}.$$
\end{proposition}
Recall the definition of  $\mathcal{C}$-subsolutions introduced by Szekleyhidi \cite{Szekelyhidi18} (see also \cite{CHZ23,HZ23}).

\begin{definition}  Let $h\in C^\infty(M)$ be a positive function and 
$$\Gamma^{h(p)} = \{ \lambda \in \Gamma~|~ f(\lambda) \geq h(p) \}, ~ p \in M. $$

\begin{enumerate}

\item   $u\in \mathcal{E}$ is said to be a $\mathcal{C}_{h} $-subsolution if the set 
$$
\left( \lambda(\omega_u (p))   + \Gamma_n \right) \cap \partial \Gamma^{h( p)} \subset \mathbb{R}^n
$$ is bounded for each $p\in M$.

\medskip

\item  $u\in  \mathcal{E}$ is said to be a $\mathcal{C}_{h, r, R}$-subsolution for some $r>0$ and $R>0$,  if 
\begin{equation}
\left( \lambda(\omega_u(p)) - r \mathbf{1} + \Gamma_n \right) \cap \partial \Gamma^{h (p)} \subset B(0, R) \subset \mathbb{R}^n
\end{equation}
for all $p\in M$, where $\mathbf{1}=(1,...,1)\in \mathbb{R}^n$. 
\end{enumerate}

\end{definition}

\begin{proposition}\label{prop3.2}
    There exists $r>0,R>0$ such that $\underline{u}$ is a $\mathcal{C}_{e^{\psi_t}+c_t,r,R}$-
    subsolution for (\ref{2.4}) for all $t\in\mathcal{T}$.
\end{proposition}
\begin{proof}
    The proof used Proposition \ref{3.1} and Lemma \ref{lem3.1}. See  \cite[Proposition 3.2]{GS24}.
\end{proof}
An argument of maximum principle will deduce the uniqueness for solutions
of equation (\ref{eq1.1}).
\begin{proposition}
    Suppose that $u\in\mathcal{E}$ solves equation (\ref{eq1.1}). Then 
    \begin{enumerate}
        \item $F[u]=\sigma e^{-\psi}$, where $\sigma$ is the sup-slope for (\ref{eq1.1});\\
        \item $\sup_{u\in\mathcal{E}}\min_M e^{-\psi}F[u]=\inf_{u\in\mathcal{E}}\max_M e^{-\psi}F[u].$ 
    \end{enumerate}
\end{proposition}
\begin{lemma}\label{lem3.5}
    Suppose $u_1,u_2\in\mathcal{E}$ are solutions of equation (\ref{eq1.1}). Then 
    $u_2-u_1=constant$.

\end{lemma}
\begin{proof}
    Let $v_t=(1-t)u_1+tu_2$. Then in a local frame, we have 
    \begin{equation*}
        \begin{aligned}
 0&=F[u_2]-F[u_1]=\int_0^1 \frac{d}{dt}F[v_t]dt \\
 &=\left(\int_{0}^{1}
    G_{i,j}[v_t]dt\right)\Big(e_{i}\bar{e}_{j}(u_2-u_1)-[e_{i},\bar{e}_{j}]^{0,1}(u_2-u_1)+e_{p}(u_2-u_1)Z_{i\bar j}^{p}+\bar{e}_{p}(u_2-u_1)\overline{Z_{i\bar j}^{p}}\Big).
    \end{aligned}
    \end{equation*}
    By the standard elliptic operator theory, $L(u_2-u_1)$ implies that $u_2-u_1$ is a constant 
    using the  strong maximum principle.
\end{proof}
Now we can prove the closedness for the continuity method (\ref{2.4}).
\begin{lemma}
    For any $k>0$, there exists $C_k>0$ such that for any $t\in\mathcal{T}$,
    the solution $\phi_t$ to equation (\ref{2.4}) satisfies $$||\phi_t||_{C^k(X)}\leq C_k.$$ 
\end{lemma}
\begin{proof}
    It follows from Proposition \ref{prop3.2} that $\underline{u}$ is a 
    $\mathcal{C}_{e^{\psi_t+c_t},r,R}$-subsolution for some fixed $0<r<R$ where
    $\psi_t$ and $c_t$ uniformly bounded with respect to $t\in\mathcal{T}$. Applying 
    a priori estimate \cite[Proposition 3.11]{HZ23}, the lemma holds.
\end{proof}
The closedness of $\mathcal{T}$ holds and  we will consequently have 
\begin{corollary}
    There exists a unique $u\in\mathcal{E}$ solving the equation 
    $$F[u]=e^{\psi+C}.$$ This proves $(4)\implies (1)$.
\end{corollary}
By a similar argument as in \cite{GS24} we have  $(1)\implies (3)\implies (2)\implies (4)$. Combining the uniqueness
result of Lemma \ref{lem3.5}, we can now complete the proof of Theorem \ref{criterion}. 

\emph{Proof of Theorem \ref{mainthm}:}
In case (1), 
observe that $w_{i\bar{j}}$, coefficience of  the form $\eta+\frac{1}{n-1}\left((\Delta^C u)\chi-\sqrt{-1}\p\bp u\right)$
can be expressed as $$w_{i\bar{j}}=P_{\chi}(\tilde{g}_{i\bar{j}})
:=\frac{1}{n-1}(tr_{\chi}\tilde{g})\chi_{i\bar{j}}-\tilde{g}_{i\bar{j}},$$
where $W_{i\bar{j}}=P_{\chi}(Z_{i\bar{j}})$ and similarly $\eta=P_{\chi}(\omega)$.
To be more precise, let $T$ be the linear invertible map defined by 
$$T(\lambda)=(T(\lambda_1,...,\lambda_n)), T(\lambda)_k=\frac{1}{n-1}\sum_{i\neq k}\lambda_i,~for~\lambda\in\mathbb{R}^n.$$
Define $$f=\left(S_n(T)\right)^{1/n},\Gamma=T^{-1}(\Gamma_n),$$
and the equation (\ref{form type eq1}) can be written as 
$$ F[u]=e^{\frac{\psi}{n}}.$$ 
In case (2), let $f=(S_k)^{1/k}$ and $\Gamma=\Gamma_k$, then the result follows.
\hfill \qedsymbol

\end{document}